\documentclass[a4paper,11pt]{amsart}
\usepackage{graphicx}
\usepackage{mathptmx}
\usepackage{amsmath}
\usepackage{amssymb}
\usepackage{enumitem}
\usepackage{xcolor}
\usepackage{xparse}
\NewDocumentCommand{\eulerian}{omm}
 {%
  \genfrac<>{0pt}{}{#2}{#3}%
  \IfValueT{#1}{_{\!#1}}%
 }

 \newmuskip\pFqmuskip

\newcommand*\pFq[6][8]{%
  \begingroup 
  \pFqmuskip=#1mu\relax
  \mathchardef\normalcomma=\mathcode`,
  \mathcode`\,=\string"8000
  \begingroup\lccode`\~=`\,
  \lowercase{\endgroup\let~}\pFqcomma
  {}_{#2}F_{#3}{\left(\genfrac..{0pt}{}{#4}{#5}\bigg|#6\right)}%
  \endgroup
}
\newcommand{\pFqcomma}{{\normalcomma}\mskip\pFqmuskip}

\newtheorem{theorem}{Theorem}

\begin{document}

\title[Some identities on $\lambda$-analogues of $r$-Stirling numbers of the second kind]{Some identities on $\lambda$-analogues of $r$-Stirling numbers of the second kind}

\author{DAE SAN KIM}
\address{Department of Mathematics, Sogang University, Seoul 121-742, Republic of Korea}
\email{dskim@sogang.ac.kr}

\author{Hyekyung Kim*}
\address{Department Of Mathematics Education, Daegu Catholic University, Gyeongsan 38430, Republic of Korea}
\email{hkkim@cu.ac.kr}

\author{Taekyun  Kim*}
\address{Department of Mathematics, Kwangwoon University, Seoul 139-701, Republic of Korea}
\email{tkkim@kw.ac.kr}

\thanks{ * are corresponding authors}
\maketitle

\begin{abstract}
Recently, the $\lambda$-analogues of $r$-Stirling numbers of the first kind were studied by Kim-Kim. The aim of this paper is to introduce the $\lambda$-analogues of $r$-Stirling numbers of the second kind and to investigate some properties, recurrence relations and certain identities on those numbers. We also introduce the $\lambda$-analogues of Whitney-type $r$-Stirling numbers of the second and derive similar results to the case of the $\lambda$-analogues of $r$-Stirling numbers of the second kind. In addition, we consider the $\lambda$-analogues of Dowling polynomials and deduce a Dobinski-like formula.
\end{abstract}

\section{Introduction}

Carlitz [3] initiated a study of the degenerate Bernoulli and Euler polynomials and numbers, which are degenerate versions of the Bernoulli and Euler polynomials and numbers.
In recent years, studying degenerate versions of special numbers and polynomials regained interests of some mathematicians. They have been explored with various tools and many fascinating results have been revealed. It is remakable that this quest for degenerate versions is not just restricted to polynomials but also extended to transcendental functions, like gamma functions. In addition, it also led to the introduction of $\lambda$-umbral calculus and $\lambda$-Sheffer sequences. \par
The degenerate Stirling numbers of the first kind and of the second kind, which are degenerate versions of the Stirling numbers of the first kind and of the second kind, appear frequently when we study degenerate versions of some special numbers and polynomials. They arise naturally when we replace the powers of $x$ by the generalized falling factorial polynomials $(x)_{k,\lambda}$ in the defining equations of the Stirling numbers of both kinds (see \eqref{1}), while the $\lambda$-analogues of Stirling numbers of the first kind and of the second kind appear when we replace the falling factorials by the generalized falling factorials.\par
The aim of this paper is to introduce the $\lambda$-analogues of $r$-Stirling numbers of the second kind and to investigate some properties, recurrence relations and certain identities on those numbers. We also introduce the $\lambda$-analogues of Whitney-type $r$-Stirling numbers of the second and derive similar results to the case of the $\lambda$-analogues of $r$-Stirling numbers of the second kind. In addition, we consider the $\lambda$-analogues of Dowling polynomials, which are a natural extension of the analogues of Whitney-type Stirling numbers of the second kind, and deduce a Dobinski-like formula.\par
The outline of this paper is as follows. In Section 1, we recall the generalized falling factorial sequence, the degenerate exponential functions, the Stirling numbers of both kinds, the $\lambda$-analogues of $r$-Stirling numbers of the first kind and the $\lambda$-analogues of unsigned $r$-Stirling numbers of the first kind. In Section 2, we introduce the $\lambda$-analogues of $r$-Stirling numbers of the second as the coefficients appearing when powers of $x+r$ are expressed in terms of the degenerate falling factorial sequence. In the special case of $r=0$, we get the $\lambda$-analogues of Stirling numbers of the second kind.  In Theorem 1, we obtain the generating function of the $\lambda$-analogues of $r$-Stirling numbers of the second. We express those numbers in terms of the forward difference operator in Theorem 2. In Theorems 3 and 6, we find an expression of the analogues of $r$-Stirling numbers of the second kind in terms of the analogues of Stirling numbers of the second kind and vice versa.
In Theorems 4 and 5, we get recurrence relations for the analogues of $r$-Stirling numbers of the second kind. In Theorem 7, we derive an identity connecting the analogues of $r$-Stirling numbers of the second kind and some values of the higher order Bernoulli polynomials. In Section 3, introduced are the $\lambda$-analogues of Whitney-type $r$-Stirling numbers of the second. In case of $r=1$, we get the $\lambda$-analogues of Whitney-type Stirling numbers of the second. The generating function of those numbers are derived in Theorem 12. In Theorem 13, we obtain an identity relating the $\lambda$-analogues of Whitney-type $r$-Stirling numbers of the second, the $\lambda$-analogues of Whitney-type Stirling numbers of the second and some values of higher order Bernoulli numbers. We introduce the $\lambda$-analogues of Dowling polynomials, and deduce the generating function and Dobinski-like formula in Theorems 10 and 11, respectively. In the rest of this section, we recall the facts that are needed in this paper.\par
Throughout this paper, let $\lambda$ be any nonzero real number. The generalized falling factorial sequence is defined by
\begin{equation}
(x)_{0,\lambda}=1,\quad (x)_{n,\lambda}=x(x-\lambda)\cdots(x-(n-1)\lambda),\quad (n\ge 1),\quad (\mathrm{see}\ [2,6,8,9,11]). \label{1}	
\end{equation}
It is known that the degenerate exponential functions are defined by
\begin{equation}
e_{\lambda}^{x}(t)=(1+\lambda t)^{\frac{x}{\lambda}}=\sum_{n=0}^{\infty}(x)_{n,\lambda}\frac{t^{n}}{n!},\quad (\mathrm{see}\ [3,7,10]). \label{2}
\end{equation}
When $x=1$, we use the notation as $e_{\lambda}(t)=e_{\lambda}^{1}(t)$. \par
For $n\ge 0$, the Stirling numbers of the first kind are defined by
\begin{equation}
(x)_{n}=\sum_{k=0}^{n}S_{1}(n,k)x^{k},\quad (\mathrm{see}\ [1,2,4,7,12,13]), \label{3}	
\end{equation}
where the falling factorial sequence is given by
\begin{equation*}
(x)_{0}=1,\quad (x)_{n}=x(x-1)\cdots(x-n+1),\quad (n\ge 1).
\end{equation*}
As the inversion formula of \eqref{3}, the Stirling numbers of the second kind are defined as
\begin{equation}
x^{n}=\sum_{k=0}^{n}{n \brace k}(x)_{k},\quad (\mathrm{see}\ [4,9,12]).\label{4}	
\end{equation}
The $\lambda$-analogues of Stirling numbers of the first kind are defined by
\begin{equation}
(x)_{n,\lambda}=\sum_{k=0}^{n}S_{1,\lambda}(n,k)x^{k},\quad (n\ge 0),\quad (\mathrm{see}\ [9]).\label{5}
\end{equation}
For $r\in\mathbb{N}\cup\{0\}$, the $\lambda$-analogues of $r$-Stirling numbers of the first kind are defined by
\begin{equation}
(x+r)_{n,\lambda}=\sum_{k=0}^{n}S_{1,\lambda}^{(r)}(n,k)x^{k},\quad (\mathrm{see}\ [9]).\label{6}
\end{equation}
Also, the $\lambda$-analogues of unsigned $r$-Stirling numbers of the first kind are given by
\begin{equation}
\langle x+r\rangle_{n,\lambda}=\sum_{k=0}^{n}{n+r \brack k+r}_{r,\lambda}x^{k},\quad (n\ge 0),\quad (\mathrm{see}\ [9]). \label{7}	
\end{equation}
Note that $\displaystyle\lim_{\lambda\rightarrow 1} {n+r \brack k+r}_{r,\lambda} = {n+r \brack k+r}_{r}\displaystyle$ are the ordinary unsigned $r$-Stirling numbers of the first kind. \par

\section{$\lambda$-analogues of $r$-Stirling numbers of the second kind}
First, we consider the $\lambda$-analogues of Stirling numbers of the second kind as the inversion formula of \eqref{5}, which are defined by
\begin{equation}
x^{n}=\sum_{k=0}^{n}{n \brace k}_{\lambda}(x)_{k,\lambda},\quad (n\ge 0). \label{8}	
\end{equation}
From \eqref{8}, we note that
\begin{equation}
\begin{aligned}
	e^{xt}=\sum_{n=0}^{\infty}\frac{t^{n}}{n!}x^{n}&=\sum_{n=0}^{\infty}\frac{t^{n}}{n!}\sum_{k=0}^{n}{n \brace k}_{\lambda}(x)_{k,\lambda} \\
	&=\sum_{k=0}^{\infty}\bigg(\sum_{n=k}^{\infty}{n \brace k}_{\lambda}\frac{t^{n}}{n!}\bigg)(x)_{k,\lambda}.
\end{aligned}	\label{9}
\end{equation}
On the other hand, by \eqref{1}, we get
\begin{align}
e^{xt}&=\big(e^{\lambda t}-1+1\big)^{\frac{x}{\lambda}}=\sum_{k=0}^{\infty}\binom{\frac{x}{\lambda}}{k}(e^{\lambda t}-1)^{k} \label{10} \\
&=\sum_{k=0}^{\infty}\lambda^{-k}\frac{1}{k!}\big(e^{\lambda t}-1\big)^{k}(x)_{k,\lambda}\nonumber.
\end{align}
From \eqref{9} and \eqref{10}, we note that
\begin{equation}
\frac{1}{\lambda^{k}}\frac{1}{k!}\big(e^{\lambda t}-1\big)^{k}=\sum_{n=k}^{\infty}{n \brace k}_{\lambda}\frac{t^{n}}{n!}.\label{11}	
\end{equation}
Note that $\displaystyle \lim_{\lambda\rightarrow 0}{n \brace k}_{\lambda}={n \brace k}\displaystyle$ are the Stirling numbers of the second kind which are defined by
\begin{equation*}
x^{n}=\sum_{k=0}^{n}{n \brace k}(x)_{k},\quad (n\ge 0).
\end{equation*}
For $r\in\mathbb{N}\cup\{0\}$, we consider the $\lambda$-analogues of $r$-Stirling numbers of the second kind as the  inversion formula of \eqref{6} which are defined by
\begin{equation}
(x+r)^{n}=\sum_{k=0}^{n}{n+r \brace k+r}_{r,\lambda}(x)_{k,\lambda},\quad (n\ge 0). \label{12}
\end{equation}
From \eqref{12}, we note that
\begin{align}
e^{(x+r)t}&=\sum_{n=0}^{\infty}(x+r)^{n}\frac{t^{n}}{n!}=\sum_{n=0}^{\infty}\bigg(\sum_{k=0}^{n}{n+r \brace k+r}_{r,\lambda}(x)_{k,\lambda}\bigg)\frac{t^{n}}{n!}\label{13}\\
&=\sum_{k=0}^{\infty}\bigg(\sum_{n=k}^{\infty}{n+r \brace k+r}_{r,\lambda}\frac{t^{n}}{n!}\bigg)(x)_{k,\lambda}.\nonumber	
\end{align}
On the other hand, by \eqref{1}, we get
\begin{align}
e^{(x+r)t}&=e^{rt}e^{xt}=e^{rt}(e^{\lambda t}-1+1)^{\frac{x}{\lambda}}\label{14} \\
&=e^{rt}\sum_{k=0}^{\infty}\binom{\frac{x}{\lambda}}{k}(e^{\lambda t}-1)^{k}=\sum_{k=0}^{\infty}\frac{1}{k}(e^{\lambda t}-1)^{k}e^{rt}\frac{1}{\lambda^{k}}(x)_{k,\lambda}. \nonumber	
\end{align}
From \eqref{13} and \eqref{14}, we note that
\begin{equation}
	\frac{1}{\lambda^{k}}\frac{1}{k!}\big(e^{\lambda t}-1\big)^{k}e^{rt}=\sum_{n=k}^{\infty}{n+r \brace k+r}_{r,\lambda}\frac{t^{n}}{n!}. \label{15}
\end{equation}
\begin{theorem}
The generating function of the $\lambda$-analogues of $r$-Stirling numbers of the second kind is given by
\begin{displaymath}
\frac{1}{\lambda^{k}}\frac{1}{k!}\big(e^{\lambda t}-1\big)^{k}e^{rt}=\sum_{n=k}^{\infty}{n+r \brace k+r}_{r,\lambda}\frac{t^{n}}{n!},\quad (k\ge 0).
\end{displaymath}
\end{theorem}
Note that $\displaystyle \lim_{\lambda\rightarrow 1}{n+r \brace k+r}_{r,\lambda}={n+r \brace k+r}_{r}\displaystyle$ are the ordinary $r$-Stirling numbers of the second kind which are defined by
\begin{equation*}
(x+r)^{n}=\sum_{k=0}^{n}{n+r \brace k+r}_{r}(x)_{k},\quad (n\ge 0),\quad (\mathrm{see}\ [\ ]).
\end{equation*}
Let $\triangle$ be the difference operator with
\begin{displaymath}
	\triangle f(x)=f(x+1)-f(x).
\end{displaymath}
Then we have
\begin{equation}
\triangle^{k}f(x)=\sum_{l=0}^{k}\binom{k}{l}(-1)^{k-l}f(x+l).\label{16}	
\end{equation}
From Theorem 1, we note that
\begin{align}
\sum_{n=k}^{\infty}{n+r \brace k+r}_{r,\lambda}\frac{t^{n}}{n!}&=\frac{1}{\lambda^{k}}\frac{1}{k!}\big(e^{\lambda t}-1\big)^{k}e^{rt} \label{17} \\
&=\frac{1}{\lambda^{k}}\frac{1}{k!}\sum_{l=0}^{k}\binom{k}{l}(-1)^{k-l}e^{(l\lambda +r)t} \nonumber \\
&=\sum_{n=0}^{\infty}\bigg(\frac{1}{\lambda^{k}}\frac{1}{k!}\sum_{l=0}^{k}\binom{k}{l}(-1)^{k-l}(l\lambda+r)^{n}\bigg)\frac{t^{n}}{n!}. \nonumber
\end{align}
Let us take $f(x)=\big(\frac{x}{\lambda}\big)^{n}$ in \eqref{8}. Then we have
\begin{align}
\triangle^{k}\bigg(\frac{r}{\lambda}\bigg)^{n}&=\sum_{l=0}^{k}\binom{k}{l}(-1)^{k-l}\bigg(l+\frac{r}{\lambda}\bigg)^{n} \label{18} \\
&=\sum_{l=0}^{k}\binom{k}{l}(-1)^{k-l}(\lambda l+r)^{n}\lambda^{-n}. \nonumber	
\end{align}
By \eqref{17} and \eqref{18}, we get
\begin{equation}
\sum_{n=k}^{\infty}{n+r \brace k+r}_{r,\lambda}\frac{t^{n}}{n!}=\sum_{n=0}^{\infty}\lambda^{n-k}\frac{1}{k!}\triangle^{k}\bigg(\frac{r}{\lambda}\bigg)^{n}\frac{t^{n}}{n!}.\label{19}	
\end{equation}
\begin{theorem}
For $k\ge 0$, we have
\begin{displaymath}
\lambda^{n-k}\frac{1}{k!}\triangle^{k}\bigg(\frac{r}{\lambda}\bigg)^{n}=\left\{\begin{array}{ccc}
{n+r \brace k+r}_{r,\lambda}, & \textrm{if $n\ge k$,} \\
0, & \textrm{if $0 \le n<k$}.
\end{array}\right.
\end{displaymath}
\end{theorem}
By Theorem 1, we get
\begin{align}
\sum_{n=k}^{\infty}{n+r \brace k+r}_{r,\lambda}\frac{t^{n}}{n!}&=\frac{1}{\lambda^{k}}\frac{1}{k!}\big(e^{\lambda t}-1\big)^{k}e^{rt} \label{20} \\
&=\sum_{l=k}^{\infty}{l \brace k}_{\lambda}\frac{t^{l}}{l!}\sum_{m=0}^{\infty}r^{m}\frac{t^{m}}{m!}\nonumber \\
&=\sum_{n=k}^{\infty}\bigg(\sum_{l=k}^{n}{l \brace k}_{\lambda}r^{n-l}\binom{n}{l}\bigg)\frac{t^{n}}{n!}.\nonumber
\end{align}
Therefore, by comparing the coefficients on both sides of \eqref{20}, we obtain the following theorem.
\begin{theorem}
For $n,k\in\mathbb{Z}$ with $n\ge k\ge 0$, we have
\begin{displaymath}
{n+r \brace k+r}_{r,\lambda}=\sum_{l=k}^{n}\binom{n}{l}{l \brace k}_{\lambda}r^{n-l}.
\end{displaymath}
\end{theorem}
For $n\ge 1$, we have
\begin{align}
(x+r)^{n+1}&=(x+r)^{n}(x+r)=\sum_{k=0}^{n}{n+r \brace k+r}_{r,\lambda}(x+r)(x)_{k,\lambda}\nonumber\\
     &=\sum_{k=0}^{n}{n+r \brace k+r}_{r,\lambda}(x-k\lambda+k\lambda+r)(x)_{k,\lambda}.   \label{21} \\
     &=\sum_{k=0}^{n}{n+r \brace k+r}_{r,\lambda}(x)_{k+1,\lambda}+\sum_{k=0}^{n}{n+r \brace k+r}_{r,\lambda}(x)_{k,\lambda}(\lambda k+r) \nonumber \\
     &=\sum_{k=1}^{n+1}{n+r \brace k-1+r}_{r,\lambda}(x)_{k,\lambda}+\sum_{k=0}^{n}(\lambda k+r){n+r \brace k+r}_{r,\lambda}(x)_{k,\lambda} \nonumber \\
     &=\sum_{k=0}^{n+1}\bigg({n+r \brace k-1+r}_{r,\lambda}+(\lambda k+r){n+r \brace k+r}_{r,\lambda}\bigg)(x)_{k,\lambda}.\nonumber
\end{align}
On the other hand, by \eqref{12}, we get
\begin{equation}
(x+r)^{n+1}=\sum_{k=0}^{n+1}{n+1+r \brace k+r}_{r,\lambda}(x)_{k,\lambda}.\label{22}	
\end{equation}
Therefore, by \eqref{21} and \eqref{22}, we obtain the following theorem.
\begin{theorem}
For $n,k\in\mathbb{Z}$ with $n\ge k\ge 1$, we have
\begin{displaymath}
	{n+1+r \brace k+r}_{r,\lambda}={n+r \brace k-1+r}_{r,\lambda}+(\lambda k+r){n+r \brace k+r}_{r,\lambda}.
\end{displaymath}
\end{theorem}
Now, we observe that
\begin{align}
&\frac{1}{\lambda^{k}}\frac{1}{k!}\big(e^{\lambda t}-1\big)^{k}\frac{1}{\lambda^{m}}\frac{1}{m!}\big(e^{\lambda t}-1\big)^{m}e^{rt}\label{23} \\
&=\frac{1}{\lambda^{k+m}}\frac{1}{(k+m)!}(e^{\lambda t}-1)^{k+m}e^{rt}\frac{(k+m)!}{k! m!} \nonumber \\
&=\binom{k+m}{k}\sum_{n=m+k}^{\infty}{n+r \brace k+m+r}_{r,\lambda}\frac{t^{n}}{n!}. \nonumber
\end{align}
On the other hand, by \eqref{15}, we get
\begin{align}
&\frac{1}{\lambda^{k}}\frac{1}{k!}\big(e^{\lambda t}-1\big)^{k}\frac{1}{\lambda^{m}}\frac{1}{m!}(e^{\lambda t}-1)^{m}e^{rt} \label{24}\\
&=\sum_{l=k}^{\infty}{l \brace k}_{\lambda}\frac{t^{l}}{l!}\sum_{j=m}^{\infty}{j+r \brace m+r}_{r,\lambda}\frac{t^{j}}{j!} \nonumber \\
&=\sum_{n=m+k}^{\infty}\sum_{l=k}^{n-m}{l \brace k}_{\lambda}{n-l+r \brace m+r}_{r,\lambda}\binom{n}{l}\frac{t^{n}}{n!}. \nonumber	
\end{align}
Therefore, by \eqref{23} and \eqref{24}, we obtain the following theorem.
\begin{theorem}
For $m,n,k\ge 0$ with $n\ge m+k$, we have
\begin{displaymath}
\binom{m+k}{k}{n+r \brace k+m+r}_{r,\lambda}=\sum_{l=k}^{n-m}\binom{n}{l}{l \brace k}_{\lambda}{n-l+r \brace m+r}_{r,\lambda}.
\end{displaymath}
\end{theorem}
From the definition of the $\lambda$-analogues of the Stirling numbers of the second kind, we have
\begin{align}
\sum_{n=k}^{\infty}{n \brace k}_{\lambda}\frac{t^{n}}{n!}&=\frac{1}{\lambda^{k}}\frac{1}{k!}\big(e^{t}-1\big)^{k}=\frac{1}{\lambda^{k}}\frac{1}{k!}\big(e^{\lambda t}-1\big)^{k} e^{rt} e^{-rt} \label{25} \\
&=\sum_{l=k}^{\infty}{l+r \brace k+r}_{r,\lambda}\frac{t^{l}}{l!}\sum_{m=0}^{\infty}(-r)^{m}\frac{t^{m}}{m!} \nonumber \\
&=\sum_{n=k}^{\infty}\bigg(\sum_{l=k}^{n}\binom{n}{l}{l+r \brace k+r}_{r,\lambda}(-1)^{n-l}r^{n-l}\bigg)\frac{t^{n}}{n!}.\nonumber	
\end{align}
Therefore, by comparing the coefficients on both sides of \eqref{25} we obtain the following theorem.
\begin{theorem}
For $n,k\ge0 $ with $n\ge k$, we have
\begin{displaymath}
{n \brace k}_{\lambda}=\sum_{l=k}^{n}\binom{n}{l}{l+r \brace k+r}_{r,\lambda}(-1)^{n-l}r^{n-l}.
\end{displaymath}
\end{theorem}
For $m\in\mathbb{N}$, the higher-order Bernoulli polynomials are defined by
\begin{equation}
\bigg(\frac{t}{e^{t}-1}\bigg)^{m}e^{xt}=\sum_{n=0}^{\infty}B_{n}^{(m)}(x)\frac{t^{n}}{n!},\quad (\mathrm{see}\ [1,3,7]). \label{26}	
\end{equation}
From \eqref{26}, we note that
\begin{align}
&\sum_{n=k}^{\infty}{n+r\brace k+r}_{r,\lambda}\frac{t^{n}}{n!}=\frac{1}{\lambda^{k}}\frac{1}{k!}\big(e^{\lambda t}-1\big)^{k}e^{rt}	\label{27} \\
&\quad =\frac{1}{t^{m}}\frac{1}{\lambda^{k+m}}
\frac{(k+m)!}{k!}\frac{1}{(k+m)!}(e^{\lambda t}-1)^{k+m}\bigg(\frac{\lambda t}{e^{\lambda t}-1}\bigg)^{m}e^{rt}\nonumber \\
&\quad =\binom{k+m}{k}\frac{m!}{t^{m}}\sum_{l=k+m}^{\infty}{l \brace k+m}_{\lambda}\frac{t^{l}}{l!}\bigg(\sum_{j=0}^{\infty}B_{j}^{(m)}\frac{r}{\lambda}\bigg)\frac{\lambda^{j}t^{j}}{j!} \nonumber \\
&\quad = \binom{k+m}{k}\sum_{l=k}^{\infty}{l+m \brace k+m}_{\lambda}\frac{m!l!}{(l+m)!}\frac{t^{l}}{l!}\sum_{j=0}^{\infty}B_{j}^{(m)}\bigg(\frac{r}{\lambda}\bigg)\frac{\lambda^{j}t^{j}}{j!} \nonumber \\
&\quad = \binom{k+m}{k}\sum_{n=k}^{\infty}\bigg(\sum_{l=k}^{n}\frac{\binom{n}{l}}{\binom{l+m}{l}}{l+m \brace k+m}_{\lambda}\lambda^{n-l}B_{n-l}^{(m)}\bigg(\frac{r}{\lambda}\bigg)\bigg)\frac{t^{n}}{n!}. \nonumber
\end{align}
Therefore, by comparing the coefficients on both sides of \eqref{27}, we obtain the following theorem.
\begin{theorem}
For $n,k\ge 0$ with $n\ge k$, we have
\begin{displaymath}
\frac{{n+r \brace k+r}_{r,\lambda}}{\binom{k+m}{k}}=\sum_{l=k}^{n}\frac{\binom{n}{l}}{\binom{l+m}{l}}{l+m \brace k+m}_{\lambda}B_{n-l}^{(m)}\bigg(\frac{r}{\lambda}\bigg)\lambda^{n-l}.
\end{displaymath}
\end{theorem}

\section{Further Remarks}
For $m,n\ge 0$, we define $\lambda$-analogues of the Whitney-type Stirling numbers of the second kind as
\begin{equation}
(mx+1)^{n}=\sum_{k=0}^{n}W_{m,\lambda}(n,k)m^{k}(x)_{k,\lambda},\quad (n\ge 0). \label{28}
\end{equation}
From \eqref{28}, we note that
\begin{equation}
\begin{aligned}
e^{(mx+1)t}&=\sum_{m=0}^{\infty}(mx+1)^{n}\frac{t^{n}}{n!}=\sum_{n=0}^{\infty}\sum_{k=0}^{n}W_{m,\lambda}(n,k)m^{k}(x)_{k,\lambda}\frac{t^{n}}{n!} \\
&=\sum_{k=0}^{\infty}\sum_{n=k}^{\infty}W_{m,\lambda}(n,k)\frac{t^{n}}{n!}m^{k}(x)_{k,\lambda}.
\end{aligned}\label{29}
\end{equation}
On the other hand, by \eqref{11}, we get
\begin{equation}
\begin{aligned}
e^{(mx+1)t}&=e^{t}\big(e^{\lambda mt}-1+1\big)^{\frac{x}{\lambda}}=e^{t}\sum_{k=0}^{\infty}\binom{\frac{x}{\lambda}}{k}(e^{\lambda mt}-1)^{k} \\
&=\sum_{k=0}^{\infty}\bigg(\frac{1}{\lambda^{k}}\frac{1}{k!}\bigg(\frac{e^{\lambda mt}-1}{m}\bigg)^{k}e^{t}\bigg) m^{k}(x)_{k,\lambda}.
\end{aligned}\label{30}
\end{equation}
Therefore, by \eqref{29} and \eqref{30}, we obtain the generating function $W_{m,\lambda}(n,k),\ (n,k\ge 0)$. \par
\begin{theorem}
For $k\ge 0$, we have
\begin{displaymath}
\frac{1}{k!}\frac{1}{\lambda^{k}}\bigg(\frac{e^{\lambda mt}-1}{m}\bigg)^{k}e^{t}=\sum_{n=k}^{\infty}W_{m,\lambda}(n,k)\frac{t^{n}}{n!}.	
\end{displaymath}
\end{theorem}
From Theorem 8, we note that
\begin{align}	
&\sum_{n=k}^{\infty}{n+1 \brace k+1}_{\lambda}\frac{t^{n}}{n!}=\frac{d}{dt}\bigg\{\sum_{n=k}^{\infty}{n+1 \brace k+1}_{\lambda}\frac{t^{n+1}}{(n+1)!}\bigg\} \label{31} \\
&=\frac{d}{dt}\bigg(\frac{1}{(k+1)!}\frac{1}{\lambda^{k+1}}(e^{\lambda t}-1)^{k+1}\bigg)=\frac{1}{k!}\frac{1}{\lambda^{k}}(e^{\lambda t}-1)^{k}e^{t}e^{(\lambda-1) t} \nonumber \\
&=\sum_{l=k}^{\infty}W_{1,\lambda}(l,k)\frac{t^{l}}{l!}\sum_{j=0}^{\infty}(\lambda -1)^{j}\frac{t^{j}}{j!} \nonumber \\
&=\sum_{n=k}^{\infty}\bigg(\sum_{l=k}^{n}W_{1,\lambda}(l,k)(\lambda -1)^{n-l}\binom{n}{l}\bigg)\frac{t^{n}}{n!}. \nonumber
\end{align}
Therefore, by comparing the coefficients on both sides of \eqref{31}, we obtain the following theorem.
\begin{theorem}
	For $n,k\ge 0$ with $n\ge k$, we have
	\begin{displaymath}
		\sum_{l=k}^{n}\binom{n}{l}W_{1,\lambda}(l,k)(\lambda -1)^{n-l}={n+1 \brace k+1}_{\lambda}.
	\end{displaymath}
\end{theorem}
Now, we consider the $\lambda$-analogues of Dowling polynomials which are defined by
\begin{equation}
d_{m,\lambda}(n,x)=\sum_{k=0}^{n}W_{m,\lambda}(n,k)x^{k},\quad (n\ge 0). \label{32}	
\end{equation}
Thus, by \eqref{32}, we get
\begin{align}
\sum_{n=0}^{\infty}d_{m,\lambda}(n,x)\frac{t^{n}}{n!}&=\sum_{n=0}^{\infty}\sum_{k=0}^{n}W_{m,\lambda}(n,k)x^{k}\frac{t^{n}}{n!}\label{33}\\
&=\sum_{k=0}^{\infty}x^{k}\sum_{n=k}^{\infty}W_{m,\lambda}(n,k)\frac{t^{n}}{n!}  \\
&=e^{t}\sum_{k=0}^{\infty}x^{k}\frac{1}{k!}\frac{1}{\lambda^{k}}\bigg(\frac{e^{\lambda mt}-1}{m}\bigg)^{k}\nonumber \\
&=e^{t}e^{x(\frac{e^{\lambda mt}-1}{\lambda m})}.\nonumber
\end{align}
\begin{theorem}
For $m\in\mathbb{N}$, we have
\begin{displaymath}
e^{t}e^{x(\frac{e^{\lambda mt}-1}{\lambda m})}=\sum_{n=0}^{\infty}d_{m,\lambda}(n,x)\frac{t^{n}}{n!}.
\end{displaymath}
\end{theorem}
When $x=1$, $d_{m,\lambda}(n,1)=d_{m,\lambda}(n)$ are called the $\lambda$-analogues of Dowling numbers. \par
We define the $\lambda$-analogues of Bell polynomials by
\begin{equation}
e^{\frac{x}{\lambda}(e^{\lambda t}-1)}=\sum_{n=0}^{\infty}\phi_{n,\lambda}(x)\frac{t^{n}}{n!}. \label{34}
\end{equation}
Thus, we easily get $\displaystyle \phi_{n,\lambda}(x)=\sum_{k=0}^{n}{n \brace k}_{\lambda}x^{k},\ (n\ge 0)\displaystyle$. When $x=1,\ \phi_{n,\lambda}(1)=\phi_{n,\lambda}$ are called the $\lambda$-analogues of Bell numbers. \par
From Theorem 10, we note that
\begin{align}
\sum_{n=0}^{\infty}d_{m,\lambda}(n,x)\frac{t^{n}}{n!}&=e^{t}e^{x(\frac{e^{\lambda mt}-1}{\lambda m})}=e^{-\frac{x}{\lambda m}}e^{t} e^{x\frac{e^{\lambda mt}}{\lambda m}}\label{37}	\\
&=e^{-\frac{x}{\lambda m}}\sum_{k=0}^{\infty}\frac{e^{\lambda mkt+t}}{\lambda^{k}m^{k}}\frac{x^{k}}{k!}\nonumber \\
&=e^{-\frac{x}{\lambda m}}\sum_{k=0}^{\infty}\frac{x^{k}}{k!m^{k}\lambda^{k}}\sum_{n=0}^{\infty}(\lambda mk+1)^{n}\frac{t^{n}}{n!} \nonumber \\
&=\sum_{n=0}^{\infty}\bigg(e^{-\frac{x}{\lambda m}}\sum_{k=0}^{\infty}\frac{x^{k}}{k!m^{k}\lambda^{k}}(\lambda mk+1)^{n}\bigg)\frac{t^{n}}{n!}.\nonumber
\end{align}
Therefore, by comparing the coefficients on both sides of \eqref{37}, we obtain the following Dobinski-like formula.
\begin{theorem}
	For $n\ge 0$, we have
	\begin{displaymath}
		d_{m,\lambda}(n,x)=e^{-\frac{x}{\lambda m}}\sum_{k=0}^{\infty}\frac{x^{k}}{k!m^{k}\lambda^{k}}(\lambda mk+1)^{n}.
	\end{displaymath}
\end{theorem}
For $r\in\mathbb{N}$, $m,n\ge 0$, we consider the $\lambda$-analogues of the Whitney- type $r$-Stirling numbers of the second kind defined by
\begin{equation}
(mx+r)^{n}=\sum_{k=0}^{n}W_{m,\lambda}^{(r)}(n,k)m^{k}(x)_{k,\lambda},\quad (n\ge 0).\label{38}
\end{equation}
From \eqref{38}, we note that
\begin{align}
e^{(mx+r)t}&=\sum_{n=0}^{\infty}(mx+r)^{n}\frac{t^{n}}{n!}=\sum_{n=0}^{\infty}\bigg(\sum_{k=0}^{n}W_{m,\lambda}^{(r)}(n,k)m^{k}(x)_{k}\bigg)\frac{t^{n}}{n!}\label{39}\\
	&=\sum_{k=0}^{\infty}\bigg(\sum_{n=k}^{\infty}W_{m,\lambda}^{(r)}(n,k)\frac{t^{n}}{n!}\bigg)m^{k}(x)_{k}.\nonumber
\end{align}
On the other hand, by \eqref{1}, we get
\begin{align}
e^{(mx+r)t}&=e^{rt}(e^{m\lambda t}-1+1)^{\frac{x}{\lambda}}=\sum_{k=0}^{\infty}\binom{\frac{x}{\lambda}}{k}(e^{m\lambda t}-1)^{k}e^{rt} \label{40} \\
&=\sum_{k=0}^{\infty}\frac{1}{k!}\frac{1}{\lambda^{k}}(e^{m\lambda t}-1)^{k}e^{rt}(x)_{k,\lambda}\nonumber \\
&=\sum_{k=0}^{\infty}\bigg(\frac{1}{k!}\frac{1}{\lambda^{k}}\bigg(\frac{e^{m\lambda t}-1}{m}\bigg)^{k}e^{rt}\bigg)m^{k}(x)_{k,\lambda}.\nonumber
\end{align}
By \eqref{39} and \eqref{40}, we get
\begin{equation}
\frac{1}{k!}\frac{1}{\lambda^{k}}\bigg(\frac{e^{\lambda mt}-1}{m}\bigg)^{k}e^{rt}=\sum_{n=k}^{\infty}W_{m,\lambda}^{(r)}(n,k)\frac{t^{n}}{n!},\label{41}
\end{equation}
where $k$ is a nonnegative integer. \par
Therefore, by \eqref{41}, we obtain the following theorem.
\begin{theorem}
For $k\ge 0$, we have
\begin{displaymath}
\frac{1}{k!}\frac{1}{\lambda^{k}}\bigg(\frac{e^{\lambda mt}-1}{m}\bigg)^{k}e^{rt}=\sum_{n=k}^{\infty}W_{m,\lambda}^{(r)}(n,k)\frac{t^{n}}{n!}.
\end{displaymath}
\end{theorem}
From \eqref{15} and Theorem 11, we have
\begin{displaymath}
	W_{1,\lambda}^{(r)}(n,k)={n+r \brace k+r}_{r,\lambda},\ W_{1,\lambda}^{(0)}(n,k)={n \brace k}_{\lambda},\ \mathrm{and}\ W_{m,\lambda}^{(1)}(n,k)=W_{m,\lambda}(n,k),\quad (n,k\ge 0).
\end{displaymath}
Now, we observe that
\begin{align}
&\sum_{n=k}^{\infty}W_{m,\lambda}^{(r)}(n,k)\frac{t^{n}}{n!}=\frac{1}{k!}\frac{1}{\lambda^{k}}\bigg(\frac{e^{\lambda mt}-1}{m}\bigg)^{k}e^{rt}\label{42}\\
&=\frac{1}{k!}\frac{1}{\lambda^{k+\alpha}}\bigg(\frac{e^{\lambda mt}-1}{m}\bigg)^{k+\alpha}e^{t}\frac{1}{t^{\alpha}}\bigg(\frac{\lambda mt}{e^{\lambda mt}-1}\bigg)^{\alpha}e^{(r-1)t}\nonumber \\
&=\frac{(k+\alpha)!}{k!}\frac{1}{t^{\alpha}}\frac{1}{(k+\alpha)!}\frac{1}{\lambda^{k+\alpha}}\bigg(\frac{e^{\lambda mt}-1}{m}\bigg)^{k+\alpha}e^{t}\sum_{j=0}^{\infty}B_{j}^{(\alpha)}\bigg(\frac{r-1}{m\lambda}\bigg)\lambda^{j}m^{j}\frac{t^{j}}{j!} \nonumber \\
&=\frac{\alpha!}{t^{\alpha}}\binom{k+\alpha}{k}\sum_{l=k+\alpha}^{\infty}W_{m,\lambda}(l,k+\alpha)\frac{t^{l}}{l!}\sum_{j=0}^{\infty}B_{j}^{(\alpha)}\bigg(\frac{r-j}{m\lambda}\bigg)\lambda^{j}m^{j}\frac{t^{j}}{j!}\nonumber \\
&=\binom{k+\alpha}{k}\bigg(\sum_{l=k}^{\infty}\frac{W_{m,\lambda}(l+\alpha,k+\alpha)}{\binom{l+\alpha}{l}}\frac{t^{l}}{l!}\bigg)\sum_{j=0}^{\infty}B_{j}^{(\alpha)}\bigg(\frac{r-j}{m\lambda}\bigg)\lambda^{j}m^{j}\frac{t^{j}}{j!}\nonumber \\
&=\binom{k+\alpha}{k}\sum_{n=k}^{\infty}\bigg(\sum_{l=k}^{n}\frac{\binom{n}{l}}{\binom{l+\alpha}{l}}W_{m,\lambda}(l+\alpha,k+\alpha)B_{n-l}^{(\alpha)}\bigg(\frac{r-n+l}{m\lambda}\bigg)\lambda^{n-l}m^{n-l}\bigg)\frac{t^{n}}{n!},\nonumber
\end{align}
where $\alpha$ is a positive integer. \par
Comparing the coefficients on both sides of \eqref{42}, we have the following theorem.
\begin{theorem}
	For $n,k\in\mathbb{N}\cup\{0\}$ and $\alpha\in\mathbb{N}$, we have
	\begin{displaymath}
	\binom{k+\alpha}{k}^{-1}W_{m,\lambda}^{(r)}(n,k)=\sum_{l=k}^{n}\frac{\binom{n}{l}}{\binom{l+\alpha}{l}}W_{m,\lambda}(l+\alpha,k+\alpha)B_{n-l}^{(\alpha)}\bigg(\frac{r-n+l}{m\lambda}\bigg)\lambda^{n-l}m^{n-l}.
	\end{displaymath}
\end{theorem}
\section{Conclusion}
In this paper, we introduced the $\lambda$-analogues of $r$-Stirling numbers of the second which appear as the coefficients when powers of $x+r$ are expressed in terms of the degenerate falling factorial sequence. We obtained some properties, recurrence relations and certain identities on such numbers. We also introduced the $\lambda$-analogues of Whitney-type $r$-Stirling numbers of the second and derived similar results to the case of the $\lambda$-analogues of $r$-Stirling numbers of the second. Furthermore, the $\lambda$-analogues of Dowling polynomials were introduced as a natural extension of the analogues of Whitney-type Stirling numbers of the second kind and a Dobinski-like formula for them was deduced.\par
It is one of our future projects to continue to study analogues of some special numbers and polynomials and to find their applications in physics, science and engineering.

\medskip

{\bf{Funding}} This work was supported by the Basic Science Research Program, the National Research Foundation of Korea,(NRF-2021R1F1A1050151).

\end{document}